\documentclass[12pt, amsfonts, epsfig]{amsart}
\usepackage{latexsym, amssymb, amsmath, epsfig, amscd, amsfonts}
 \newlength{\baseunit}               % the basic unit length
 \newcount{\numlines}                % depth of picture (in number of lines)
\newcommand{\N}{\mathbb{N}}

\newcommand{\Q}{\mathbb{Q}}
\newcommand{\R}{\mathbb{R}}
\newcommand{\C}{\mathbb{C}}

\newcommand{\proj}{\mathbb P}

        \newfont{\hollow}{msbm10 scaled\magstep1}
        
        \newfont{\Bfmit}{eufm10 scaled\magstep1}

\newcommand{\cO}{{\mathcal O}}

\newcommand{\SL}{\operatorname{SL}}

\newcommand{\id}{\operatorname{id}}

\newcommand{\SU}{\operatorname{SU}}
\newcommand{\SO}{\operatorname{SO}}

\newcommand{\Proj}{\operatorname{Proj}}

\newtheorem{thm}{Theorem}[subsection]

\newtheorem{lem}[thm]{Lemma}
\newtheorem{cor}[thm]{Corollary}

\theoremstyle{definition}
\newtheorem{defn}[thm]{Definition}

\newtheorem{rem}[thm]{Remark}           %\renewcommand{\theremark}{}

\theoremstyle{remark}

\newcommand{\lremind}[1]{{}}
\newcommand{\bremind}[1]{{}}

\newcommand{\cut}[1]{}
\begin{document}
\pagestyle{plain} \title{{ \large{Quotients by Reductive Group,
Borel Subgroup, Unipotent Group and Maximal Torus} } }
\author{Yi Hu }

\address{Department of Mathematics, University of Arizona, Tucson,
AZ 86721, USA}

\email{yhu@math.arizona.edu}

\address{Center for Combinatorics, LPMC, Nankai University, Tianjin 300071,
China}
%IMS, Chinese University of Hong Kong, Hong Kong, China}
%\email{yhu@ims.cuhk.edu.hk}

\begin{abstract} Consider an algebraic action of a connected complex
reductive algebraic group on a  complex polarized projective
variety. In this paper, we first introduce the nilpotent quotient,
the quotient of the polarized projective variety by a maximal
unipotent subgroup. Then, we introduce and investigate three
induced actions: one by the reductive group, one by a Borel
subgroup, and one by a maximal torus, respectively. Our main
result is  that there are natural correspondences among quotients
of these three actions. In the end, we mention a possible
application to the moduli spaces of parabolic bundles over
algebraic curves for further research.
\end{abstract}

\maketitle

\hskip 2cm {\scriptsize {\it ---  Dedicated to Robert MacPherson
on the occasion of his 60th birthday}}

\vskip 1cm

\section{Introduction and Statements of Results}

Let $G \times X \rightarrow X$ be an algebraic action of a
connected complex reductive algebraic group $G$ on an arbitrary
complex projective variety $X$. Let $L$ be a very ample line
bundle over $X$. We assume that  $L$ admits a
$G$-linearization\footnote{This is automatically satisfied if $X$
is normal.}.

Under these assumptions, we will introduce three other actions and
study relations among their quotients.

To this end, we fix a Borel subgroup $B$ of $G$, the unipotent
radical $U$ of $B$,  and a maximal torus $H$ of $G$ such that $B =
HU$. We also fix a compact form $K$ of $G$ such that $T= K \cap H$
is the compact torus of $H$.  Let ${\frak t}^*$ be the linear dual
of ${\frak t}= {\rm Lie} (T)$ and ${\frak t}^*_+$ be the closed
Weyl chamber in ${\frak t}^*$ which is positive with respect to
$B$.

\subsection{\it Nilpotent quotient.} Our first theorem is that
there exists a canonically defined quotient by the unipotent group
$U$.

\begin{thm}
\label{unipotentQuotient}
 There is a uniquely defined Zariski
open subset $X^{ss}_U(L)$ of $X$, which solely depends on $L$ but
not on the linearization of $L$,  such that $X^{ss}_U(L^n) =
X^{ss}_U(L)$ for all $n
>0$ and the quotient $X^{ss}_U(L)
\rightarrow X^{ss}_U(L)/\!/U$ exists. Furthermore, on the quotient
variety $X^{ss}_U(L)/\!/U$ the maximal torus $H$ naturally acts.
\end{thm}

\subsection{\it Three actions.}  The three actions that we mentioned earlier are:
\begin{enumerate}
\item   the induced torus action
$$ H \times X^{ss}_U(L)/\!/U \rightarrow X^{ss}_U(L)/\!/U;$$
\item the induced diagonal action
$$G\times (X \times G/B) \rightarrow (X \times G/B);$$
\item  and the induced Borel subgroup action
$$B \times X \rightarrow X.$$
 \end{enumerate}

 To explain the
natural correspondences among  quotients of these three actions,
we will start with a parameter space, a rational polytope
$\Delta$,  for these quotients.

\subsection{\it A parameter space for quotients.}
Let $\Lambda = \hbox{Hom}(T,\hbox{U}(1))$ be the weight lattice of
$T$ and $\Lambda_+ = \Lambda \cap {\frak t}^*_+$. (Here, as usual,
we will identify $\Lambda$ with a subgroup of ${\frak t}^*$ by
identifying the weight $\lambda$ with the functional ${\rm d}
\lambda/(2\pi i)$.) Let $({\frak t}^*_+)_{\rm reg}$ be the set of
regular points of ${\frak t}^*_+$, that is, the set of points
outside the walls of Weyl chamber.

Choose a $K$-invariant Hermitian form on the space of global
sections of $L$ and let $\Phi: X \rightarrow {\frak k}^*$ be the
associated moment map\footnote{$\Phi$ is the restriction of the
corresponding moment map on the projective space
$\proj(H^0(X,L)^*)$. Hence it makes sense even if $X$ is
singular.}
 where ${\frak k}^*$ is the linear dual of ${\frak k}= {\rm Lie}(K)$.
Let $\Delta= \Phi (X) \cap {\frak t}^*_+$. This is a rational
convex polytope (see Mumford's Appendix to \cite{Ness}. See also
\cite{Kir} for the case of symplectic manifolds).
%Brion, in  \cite{Brion},  described the rational points
%of $\Delta$ in terms of $G$-submodules of the section ring
%$\bigoplus_n H^0(X, L^n)$.
Set
$$\Delta_{\rm reg}= \Delta \cap ({\frak t}^*_+)_{\rm reg}.$$

For any rational point ${\chi \over n} \in \Delta_{\rm reg}$ with
$\chi \in \Lambda_+$ and $n \in \N$, we will associate a quotient
for each of the above three actions as follows.

\subsection{\it The first group action.}  For a sufficiently divisible
positive integer $n$, $L^n$ descends to a very ample line bundle
$\cO_{X^{ss}_U/\!/U}(n)$ over $X^{ss}_U/\!/U$ on which $H$ acts
linearly. Let $\cO_{X^{ss}_U/\!/U}(n, \chi)$ be the linearization
on $\cO_{X^{ss}_U/\!/U}(n)$  shifted by the character $-\chi$ (see
\S \ref{shiftlinearization} for the precise definition). We will
denote the corresponding locus of $H$-semistable points by
$(X^{ss}_U(L)/\!/U)_H^{ss}(n,\chi)$. This leads to our first (GIT)
quotient
$$(X^{ss}_U(L)/\!/U)_H^{ss}(n,\chi) \rightarrow (X^{ss}_U(L)/\!/U)_H^{ss}(n,\chi)/\!/H.$$

\subsection{\it The second group action.}
\label{2ndAction} Let $\C_{-\chi}$ be the one-dimensional
$B$-module with character $-\chi$ and $L_\chi = G \times_B
\C_{-\chi}$ be the corresponding linearized ample line bundle over
$G/B$. Then $L^n \otimes L_\chi$ becomes a linearized ample line
bundle over $X \times G/B$. This gives rise to our second (GIT)
quotient
$$(X \times G/B)_G^{ss}(L^n \otimes L_\chi) \rightarrow (X \times G/B)_G^{ss}(L^n
\otimes L_\chi)/\!/G.$$
% We will denote this quotient by $Q_{n, \chi}$.

\subsection{\it The third group action.}
\label{3rdAction}  Define a morphism
$$\iota_B: X \rightarrow X \times G/B, \;\;\;
\iota_B: x \rightarrow (x, [B]), \; \forall x \in X.$$  Set
$X_B^{ss}(n, \chi) = i_B^{-1}(i_B(X) \cap (X \times G/B)^{ss}(L^n
\otimes L_\chi)) \subset X.$ Then, we will show (Theorem
\ref{B-categoricalquotient}) that $X_B^{ss}(n, \chi)$ is
$B$-invariant, Zariski open, and admits a categorical quotient by
the Borel subgroup $B$
$$X_B^{ss}(n, \chi) \rightarrow X_B^{ss}(n, \chi)/\!/B.$$

\subsection{\it The correspondences.}
Here comes our second main theorem.

\begin{thm}
\label{mainthm}
 For every  rational point ${\chi \over
n}$ in $\Delta_{\rm reg}$\footnote{Here we indicate that the
rational points of $\Delta$ that are {\it not regular} will have
to be treated separately as they are related to homogeneous spaces
$G/P$ where $P$ is a parabolic subgroup strictly containing $B$.
See \S \ref{G/P} for details.}, there exists a quotient variety
for each of the three actions listed as follows:
$$(X^{ss}_U(L)/\!/U)_H^{ss}(n,\chi)/\!/H,$$
$$(X \times G/B)_G^{ss}(L^n \otimes L_\chi)/\!/G,$$
$$X_B^{ss}(n, \chi)/\!/B.$$
Moreover, these quotients as projective varieties are all
naturally isomorphic to each other.
\end{thm}

Here we mention that these correspondences hold over an arbitrary
ground field of characteristic zero. Working over the field of
complex numbers is only for the interpretation of the polytope
$\Delta$ in terms of moment map. But the use of moment map,
although convenient and adding some symplectic flavors to the
work, can be completely avoided. For example, to avoid the use of
moment map in this introduction, we could have simply used the
moment-map-free descriptions of the rational points of $\Delta$ as
in Equation (\ref{momentmapfree}) of \S \ref{revisited}.

%The detailed presentations of these correspondences will occupy
%the rest of the paper.

\section{The quotient of $X$ by the unipotent subgroup $U$}
\label{byuni}

\subsection{\it $U$-invariants of the section algebra.}
Since the line bundle $L$ is  $G$-linearized, we have that
$H^0(X,L)$ is a $G$-module, so is its linear dual $V=H^0(X,L)^*$.
Since $L$ is very ample, a choice of a basis of $H^0(X,L)^*$ will
equivariantly embed $X$ into the projective space $\proj (V)$.

Consider the $\N$-graded section algebra
$$R = \bigoplus_{d \ge 0}  R_d = \bigoplus_{d \ge 0} H^0(X,L^d)$$
on which $G$, hence $U$, acts linearly. Let
$$R^U = \bigoplus_{d \ge 0} R_d^U = \bigoplus_{d \ge 0} H^0(X,L^d)^U$$
be the subalgebra of $U$-invariant sections with the induced
grading by $\N$. Then, $R^U$ is finitely generated. To see this,
let
$$S= \bigoplus_{d \ge 0}  S_d = \bigoplus_{d \ge 0} H^0(\proj (V),\cO_{\proj (V)} (d))$$
be the polynomial algebra.
Let %$I_X$ be the homogeneous ideal of $X \subset \proj (V)$ and let
$$\pi: S \rightarrow R$$
be the restriction homomorphism. Then $R$ is a finite $S$-module.
Hence $R$ is finitely generated. Then by \cite{Grosshans}, $R^U$
is  finitely generated as well.
%As the restriction map from
%sections of $\cO_{\proj (V)} (d)$ on the whole projective space to
%sections on $X$, $\pi$ may fail to be surjective (this is related
%to projective normality of $X$). But this map is surjective in
%large degrees, and hence $R$ is always a finite $S$-module. Since
%it is known that $S^U$ is finitely generated, this
%implies that finite generation of $R^U$. %then $\hbox{ker} \;
%\pi$ is isomorphic to $I_X$ and $R$ is isomorphic to $S/I_X$.
%Moreover, $\pi$ is a $U$-homomorphism and therefore $R^U$ is
%isomorphic to $(S/I_X)^U$. Hence $R^U$ is also isomorphic to
%$S^U/I_X^U$ (this may fail over positive characteristic).
%hence $R^U \cong S^U/I_X^U$ is finitely generated.
(I thank Michel Brion for  pointing out the reference
\cite{Grosshans}.)

\subsection{\it The unipotent quotient.}
\label{unipotentquotient}
 Here comes  our main
definition.

\begin{defn}
The quotient of $X$ by the unipotent group $U$ with respect to the
linearization $L$ is defined to be $\Proj(R^U).$
\end{defn}

\subsection{\it Proof of Theorem \ref{unipotentQuotient}.}  Set
$$X^{ss}_U(L) = \{x \in X | \exists \; d > 0, s \in R_d^U, s(x) \ne 0 \}.$$
Then there is a quotient map
$$X^{ss}_U(L) \rightarrow \Proj(R^U),$$
locally induced from the inclusions $R^U \subset R$ over affine
patches $s(x) \ne 0$. Hence we will also denote $\Proj(R^U)$ by
$X^{ss}_U(L)/\!/U$.

The equality $X^{ss}_U(L^n)= X^{ss}_U(L)$ for all $n >0$ follows
immediately from the definition.

To show that $X^{ss}_U(L)$ is independent of the linearization of
$L$, note that if we change the linearization of $L$, then the
corresponding $G$-linear actions on the section algebra $R$ only
differ by  shifting a character of $G$. Since the character is
trivial when restricted to $U$, the action of $U$ on $R$ remains
unchanged. Therefore $R^U$, hence also $X^{ss}_U(L)$, only depends
on the underlying line bundle $L$ but not the linearization.

The maximal torus $H$ obviously acts linearly on $R^U$. Hence it
acts on the quotient $\Proj(R^U)= X^{ss}_U(L)/\!/U$.

This complete the proof of Theorem \ref{unipotentQuotient}.

\begin{rem}
\label{distintqUq} The unipotent quotient $X^{ss}_U(L)/\!/U$ in
general does depend on the choice of the underlying line bundle
$L$. We will justify this assertion in \S \ref{moreonexample}.
\end{rem}

\begin{rem}
It would be nice if $X^{ss}_U(L) \rightarrow \Proj(R^U)$ is a
categorical quotient (see Theorem 0.5 of \cite{GIT} for the
definition of categorical quotient). We are not able to prove this
although we believe this is true. It is worth mentioning that the
rest of the quotients considered in this paper are all
categorical. We also expect that the unipotent quotient
$\Proj(R^U)$ should admit other interpretations and bear
interesting applications.
\end{rem}

\subsection{\it A theorem of Guillemin and Sjamaar.}
By Guillemin and Sjamaar (\cite{GSj}), the unipotently semistable
locus $X^{ss}_U(L)$ admits the following description:

\begin{thm}
\label{GSj} {\rm (Theorem 4.2, \cite{GSj})}
$$X^{ss}_U(L) = \{x \in
X | \Phi (\overline{B \cdot x}) \cap {\frak t}^*_+ \ne \emptyset
\}.$$
\end{thm}

We will not use this result, except in Example \ref{1stexmp}
below.

\subsection{\it An example.}
\label{1stexmp}
 Consider the diagonal action of $G=\SL(2,\C)$ on
$(\proj^1)^n$. Let $B$ be the subgroup of upper triangular
matrixes and $U$ be the unipotent radical.

 We represent a point of
$\proj^1$ by ${\scriptsize \left[
\begin{array}{c}
          a \\
         b
   \end{array}
 \right]}.
$ Then $U$ fixes the point ${\scriptsize \left[
\begin{array}{c}
          1 \\
         0
   \end{array}
 \right]}
$ and
$$\proj^1 \setminus {\scriptsize \left[
\begin{array}{c}
          1 \\
         0
   \end{array}
 \right]}
$$
 is a single $U$-orbit on which $U$ acts freely.

We will identify
 the linear dual of the Lie algebra of
 $\SU (2)= \SO(3) $ with $\R^3$. Using a coadjoint orbit,
 we will realize $\proj^1$ as the unit sphere $S^2$ in $\R^3$.
 Under this identification, the moment map
 is simply the inclusion:  $S^2 \subset \R^3$.
 Let $p = S^2 \cap {\frak t}^*_+$.
Then $S^2$ is the coajoint orbit through $-p$. Under the
identification $G/B = K/T =  S^2$ (cf. the paragraph around
Equation (1) of \cite{GSj}), we have $[B]=[T]=-p$. Hence $-p$ is
fixed by the action of $B$. It follows that $-p$
 is  ${\scriptsize \left[
\begin{array}{c}
          1 \\
         0
   \end{array}
 \right]}.$
 Then $p$, as the only other fixed point of the maximal torus, must be
${\scriptsize \left[
\begin{array}{c}
          0 \\
         1
   \end{array}
 \right]}$.

 Let $d_i $ ($1 \le i \le n$) be some positive integers and let
 $L$ be the ample line bundle $\otimes_i \cO(d_i)$ over $(\proj^1)^n$. Then the induced moment
 map $\Phi$ is simply $\sum_i d_i \Phi_i$ where $\Phi_i$ is the
 following
composition map: the projection of
  $(\proj^1)^n$ to the ith factor followed by the inclusion $S^2 \subset \R^3$.

 Assume that $d_n$ is sufficiently large relative to other
$d_i$ ($1 \le i \le n-1$). Then by applying Guillemin-Sjamaar's
Theorem (Theorem \ref{GSj}), it is straightforward to check that
 $$ {\scriptsize \left[
\begin{array}{ccccccccc}
          a_1  &  \cdots & a_{n-1} & a_n  \\
          b_1  & \cdots & b_{n-1} & b_n
   \end{array}
 \right]} \in X^{ss}_U(L) \Longleftrightarrow b_n \ne 0.$$

We will represent an arbitrary point of $X^{ss}_U(L)$ by
$${\scriptsize
  \left[
\begin{array}{ccccccccc}
          \cdots  & a_{i_1} & \cdots & a_{i_r}  & \cdots & a_n \\
          \cdots &  1      & \cdots & 1 &  \cdots & 1
   \end{array}
 \right]}
$$
where the dotted columns are all $ {\scriptsize \left[
\begin{array}{c}
          1 \\
         0
   \end{array}
 \right]}.$
 Such a representation is obviously unique. Now define a morphism
$$\phi: X^{ss}_U(L) \rightarrow (\proj^1)^{n-1}$$
by
$${\scriptsize \left[
\begin{array}{ccccccccc}
          \cdots  & a_{i_1} & \cdots & a_{i_r}  & \cdots & a_n \\
          \cdots &  1      & \cdots & 1 &  \cdots & 1
   \end{array}
 \right] } \to
{\scriptsize \left[
\begin{array}{ccccccccc}
          \cdots  & a_{i_1}-a_{i_2} & \cdots & a_{i_r} - a_n  & \cdots  \\
          \cdots &  1      & \cdots & 1 &  \cdots
   \end{array}
 \right]},
 $$
 where the dotted columns stay the same, that is, are all ${\scriptsize \left[
\begin{array}{c}
          1 \\
         0
   \end{array}
 \right]}.$ (The column ${\scriptsize \left[
\begin{array}{c}
          a_n \\
         1
   \end{array}
 \right]}$ is deleted by the map $\phi$.)
 Then one checks easily that $\phi$ is surjective and
 $U$-equivariant where $U$ acts on the image $(\proj^1)^{n-1}$
 trivially.

  To see that $\phi$ sends distinct
 orbits to distinct points, suppose that we have
 $$\phi({\scriptsize
  \left[
\begin{array}{ccccccccc}
          \cdots  & a_{i_1} & \cdots & a_{i_r}  & \cdots & a_n \\
          \cdots &  1      & \cdots & 1 &  \cdots & 1
   \end{array}
 \right]}) = \phi({\scriptsize
  \left[
\begin{array}{ccccccccc}
          \cdots  & b_{i_1} & \cdots & b_{i_{r'}}  & \cdots & b_n \\
          \cdots &  1      & \cdots & 1 &  \cdots & 1
   \end{array}
 \right]}).$$
 Then we must have $r=r'$ and
 $$a_{i_j}-a_{i_{j+1}} = b_{i_j}-b_{i_{j+1}}, \;\;\; 1 \le j \le r$$
 where we set $a_{i_{r+1}} = a_n$ and $b_{i_{r+1}} = b_n$.
 This implies that
$$b_{i_j}-a_{i_j} = b_{i_{j+1}}-a_{i_{j+1}}, \;\;\; 1 \le j \le r.$$
Set $x=b_{i_j}-a_{i_j}$ for any $1 \le j \le r$, and let
$$u ={\scriptsize \left(
\begin{array}{cc}
          1  &  x  \\
          0  & 1
   \end{array}
 \right)}.$$
 Then
$${\scriptsize
  \left[
\begin{array}{ccccccccc}
          \cdots  & b_{i_1} & \cdots & b_{i_r}  & \cdots & b_n \\
          \cdots &  1      & \cdots & 1 &  \cdots & 1
   \end{array}
 \right]} = u \cdot {\scriptsize
  \left[
\begin{array}{ccccccccc}
          \cdots  & a_{i_1} & \cdots & a_{i_r}  & \cdots & a_n \\
          \cdots &  1      & \cdots & 1 &  \cdots & 1
   \end{array}
 \right]}.$$
 That is, $\phi: (\proj^1)^n \setminus {\scriptsize \left[
\begin{array}{ccccccccc}
          1  &  \cdots & 1  \\
          0  & \cdots & 0
   \end{array}
 \right]} \rightarrow
 (\proj^1)^{n-1}$ is a quotient map and $(\proj^1)^{n-1}$
 parameterizes the $U$-orbits on $X^{ss}_U(L).$

 %From here, one can check further that $\phi$ is universal in the sense that
 %any $U$-invariant morphism
 %$$\gamma: (\proj^1)^n \setminus {\scriptsize\left[
%\begin{array}{ccccccccc}
%          1  &  \cdots & 1  \\
%          0  & \cdots & 0
%   \end{array}
% \right]}  \rightarrow Z$$
% where $U$ acts trivially on $Z$ factors through $$\phi: (\proj^1)^n \setminus {\scriptsize \left[
%\begin{array}{ccccccccc}
%          1  &  \cdots & 1  \\
%          0  & \cdots & 0
%   \end{array}
% \right] }\rightarrow (\proj^1)^{n-1},$$ that is, the quotient
% $\phi$ is categorical (consult Definition \ref{categorical}).

Similarly, for every $ 1 \le i \le n$, by assuming that $d_i$ is
sufficiently large relative to the rest, we will get
$$ {\scriptsize \left[
\begin{array}{ccccccccc}
          a_1  &  \cdots & a_{n-1} & a_n  \\
          b_1  & \cdots & b_{n-1} & b_n
   \end{array}
 \right]} \in X^{ss}_U(L) \Longleftrightarrow b_i \ne 0$$ and its
 quotient by $U$ can also be identified with $(\proj^1)^{n-1}$.

%\begin{exmp}
%$\SL(2,\C)$ acts on $\proj(\C^2 \oplus \C^2)$.
%\end{exmp}

\section{Quotients of $X^{ss}_U(L)/\!/U$  by $H$}
\label{torusaction}

The maximal torus $H$ acts on $X^{ss}_U(L)/\!/U = \Proj(R^U)$ via
the induced linear action on $R^U$. We now study the $H$-quotients
on $X^{ss}_U(L)/\!/U$.

\subsection{\it $R^U$ as $H$-modules.}
The algebra $R$ is also $(\N \times \Lambda_+)$-graded:
$$R = \bigoplus_{d, \tau}  R_{d, \tau}$$
where $R_{d, \tau}$ is the isotypical $G$-submodule of $R_d$ of
highest weight $\tau$.

The algebra of $U$-invariant, $R^U$, inherits an $\N \times
\Lambda_+$-grading
$$R^U = \bigoplus_{d \in \N, \tau \in \Lambda_+}  R_{d, \tau}^U.$$
The maximal torus $H$ acts on $R^U$, having  $R_{d, \tau}^U$ as
the weightspace with weight $\tau$.

\subsection{\it The parameter space $\Delta$, revisited.}
\label{revisited} The rational points of the polytope $\Delta$ can
be determined purely algebraically as follows (see Mumford's
appendix to \cite{Ness} and Brion's paper \cite{Brion}): For any
$\tau \in \Lambda_+$, $d \in \N$,
$${\tau \over d} \in \Delta \Longleftrightarrow R_{d, \tau} \ne 0 $$
Alternatively, let $\Delta_\Q$ denote the set of rational points
in $\Delta$, then we have
\begin{equation}
\label{momentmapfree} \Delta_\Q = \{\;{\tau \over d}\; | \;R_{d,
\tau} \ne 0 \; \}.
\end{equation}

\subsection{\it Shifting the linearization.}
\label{shiftlinearization} For any ${\chi \over n} \in \Delta_{\rm
reg}$, we can shift the $H$-action on $L^n$ by the character
$-\chi$. In terms of the action on the section algebra of $L^n$,
the new linear action of $H$ is defined as follows: $H$ acts on
the weighspace $R_{nd, \tau}$ with weight $\tau-d\chi$ for all $d$
and $\tau$. We will denote the new $H$-linearized line bundle by
$L^n[\chi]$\footnote{A remark on notations: the character between
the brackets, e.g., $L^n[\chi]$, always indicates a shifting of a
linear action. However, $L_\chi$ is the line bundle over the flag
variety  $G/B$ and has nothing to do with shifting of
linearization.}. It is worth mentioning that the shifting does not
affect the $U$-action on the section algebra of $L^n$ because any
character is trivial when restricted to $U$. But it obviously does
affect the $H$-action on the section algebra of $L^n$ and hence
also the $B$-action on the section algebra of $L^n$.

For $L^n$ with $n$ sufficiently divisible, it descends to a very
ample line bundle $\cO_{X^{ss}_U/\!/U}(n)$ over $ X^{ss}_U/\!/U$
with an induced linear action by the maximal torus $H$. Likewise,
the linearized line bundle $L^n[\chi]$ also descends to a
$H$-linearized line bundle over $ X^{ss}_U/\!/U$, which we will
denote by $\cO_{X^{ss}_U/\!/U}(n, \chi)$. In terms of linear
actions on the section algebra, $H$ acts on $R^U_{dn, \tau}$ with
weight $\tau - d \chi$.

Denote the section algebra of $L^n$ by
$$R_{(n)}= \bigoplus_{d \ge 0}  R_{nd} = \bigoplus_{d \ge 0} H^0(X,L^{nd}).$$
Then we will use $R_{(n)}^{H[\chi]}$ and $R_{(n)}^{B[\chi]}$ to
denote the $H$ and $B$-invariants of $R_{(n)}$ under the
$(-\chi)$-shifting, respectively.

\subsection{\it $H$-Quotients of $X^{ss}_U/\!/U$.}

\begin{thm}
\label{Hquotient}  With respect to the linearized ample line
bundle $\cO_{X^{ss}_U/\!/U}(n, \chi)$, the GIT quotients
$(X^{ss}_U(L)/\!/U)_H^{ss}(n,\chi)/\!/H$ is
$$\Proj ((R_{(n)}^U)^{H[\chi]}) = \Proj (\bigoplus_d R^U_{nd, d\chi}).$$
\end{thm}
\proof By the (original) induced $H$-action on
$\cO_{X^{ss}_U/\!/U}(n)$ , we have that $R_{(n)}^U$ decomposes
into a direct sum of $H$-submodules
$$R_{(n)}^U = \bigoplus_{d, \tau}  R^U_{nd, \tau}.$$
Under the $(-\chi)$-shifted linear action, $H$ acts on the
weighspace $R^U_{nd, \tau}$ with weight $\tau - d\chi$, hence we
obtain
$$(R_{(n)}^U)^{H[\chi]} = \bigoplus_d R^U_{nd, d\chi}.$$
The statement of the theorem then follows readily.
%$$R_{(n)}^{B[\chi]} = (R_{(n)}^U)^{H[\chi]} = \bigoplus_d R^U_{nd, d\chi}$$
\endproof

\begin{rem} For sufficiently divisible $n$,  $n \Phi (X) \cap {\frak
t}^*_+$ is an integral polytope. Hence by Atiyah's version of the
Atiyah-Guillemin-Sternberg convexity theorem (\cite{Atiyah}), we
expect that under a suitable $H$-equivariant projective embedding
of $X^{ss}_U/\!/U$, the image of the induced $H$-moment map on
$X^{ss}_U/\!/U$ should precisely be $n \Phi (X) \cap {\frak
t}^*_+$.
\end{rem}

\section{Quotients of $X \times G/B$ by $G$}
\label{byG}

In this section, we will basically recollect some known results
that will be useful for our purposes.

\subsection{\it Moment maps on $G/B$ and coadjoint orbits.}
Recall (see, e.g.,  \cite{GSj}) that for any $\chi \in \Lambda
\cap ({\frak t}^*_+)_{\rm reg}$, let $\C_{-\chi}$ be the
one-dimensional $B$-module with character $-\chi$, then $L_\chi =
G \times_B \C_{-\chi}$ is a $G$-linearized ample line bundle over
$G/B$. The curvature from $\omega_\chi$ (with respect to the
$G$-invariant Hermitian metric on $L_\chi$ defined by the usual
norm on $\C$) is K\"ahler.

For ${\chi \over n} \in \Delta_{\rm reg}$, we will consider the
K\"ahler manifold
$$(G/B, \omega_{\chi \over n})$$ where $\omega_{\chi \over n}= {1
\over n} \omega_\chi$. The induced moment map is found by
composing the maps
$$ G/B \rightarrow K/T \rightarrow {\frak t}^*$$ where the first
map is the inverse of the diffeomorphism $K/T \rightarrow G/B$
induced by the inclusion and the second map is defined by $$[kT]
\rightarrow k \cdot (-{\chi \over n}).$$ In fact, this gives rise
to a symplectomorphism from $(G/B, \omega_{\chi \over n})$ to the
coadjoint orbit through $-{\chi \over n}$, ${\sl O}_{-{\chi \over
n}}$.

\subsection{\it The shifting trick and GIT quotients.}
\label{shiftingtrick} Let $\bar{\sl O}_{-{\chi \over n}}$ denote
the symplectic manifold obtained from the symplectic manifold
${\sl O}_{-{\chi \over n}}$ by replacing its symplectic form
$\omega_{\chi \over n}$ by $-\omega_{\chi \over n}$. Then the
product symplectic manifold $X \times \bar{\sl O}_{-{\chi \over
n}}$ admits a moment map $$ \widetilde{\Phi}: X \times \bar{\sl
O}_{-{\chi \over n}} \rightarrow {\frak k}^*$$ defined by the
formula
$$
\widetilde{\Phi} (x,q) = \Phi (x)-q.
$$
Now the set $\widetilde{\Phi}^{-1}(0)$ becomes identified with
the set $\Phi^{-1}({\sl O}_{-{\chi \over n}})$ and we obtain the
following identifications
$$
\widetilde{\Phi}^{-1}(0)/K = \Phi^{-1}({\sl O}_{-{\chi \over
n}})/K = \Phi^{-1}(-{\chi \over n})/K_{-{\chi \over n}}
$$
where $K_{-{\chi \over n}}$ is the isotropy subgroup of $K$ at
$-{\chi \over n}$. The above is the so-called shifting trick
(between  the symplectic reduction at a general coadjoint orbit
${\sl O}_{-{\chi \over n}}$ and  the symplectic reduction at the
origin).

The following theorem was formulated in Theorem 2.2.4 of
\cite{DH}.  It basically  follows from Mumford's Appendix to
\cite{Ness} and Theorem 8.3 of \cite{GIT}.

\begin{thm}  Let $(X \times G/B)_G^{ss}(L^n
\otimes L_\chi)$ be the semistable locus in $X \times G/B$ with
respect to the $G$-linearized line bundle $L^n \otimes L_\chi$.
Then we have a natural homeomorphism from  $\Phi^{-1}({\sl
O}_{-{\chi \over n}})/K$ to $(X\times G/B)^{ss} (L({\chi\over
n}))/\!/G$.
\end{thm}

\begin{rem}
\label{orbifolds} It follows from Theorem 8.3 of \cite{GIT} that
when ${\chi \over n} \in \Delta_{\rm reg}$ is a regular value of
the moment map $\Phi$, $(X \times G/B)_G^{ss}(L^n \otimes L_\chi)$
consists of only stable points, hence the quotient $(X\times
G/B)^{ss} (L({\chi\over n}))/\!/G$ has at worst finite quotient
singularities in this case.
\end{rem}

\section{Quotients of $X$ by  $B$}
\label{byB}

\subsection{\it The Zariski open subset $X_B^{ss}(n, \chi)$.}
\label{B-opensubset}
 Recall from \S \ref{3rdAction} that we have the morphism
$$\iota_B: X \rightarrow X \times G/B, \;\;
\iota_B: x \rightarrow (x, [B]), \; \forall x \in X.$$ This embeds
$X$ into $X \times G/B$ as the fiber over the base point $[B] \in
G/B$. (It is easy to see the morphism $\iota_B$ induces a
bijection between the set of $B$-orbits on $X$  and the set of
$G$-orbits on $X \times G/B$. Hence it is simply natural to expect
$B$-quotients on $X$ should correspond to $G$-quotients on $X
\times G/B$.)

As before, we have ${\chi \over n} \in \Delta_{\rm reg}$ with
$\chi \in \Lambda_+$ and $n \in \N$. Set
$$X_B^{ss}(n, \chi)= \{x \in X | (x, [B]) \in (X \times G/B)_G^{ss}(L^n \otimes L_\chi)\}.$$
That is, $$X_B^{ss}(n, \chi) = i_B^{-1}(i_B(X) \cap (X \times
G/B)_G^{ss}(L^n \otimes L_\chi)).$$ Clearly, $X_B^{ss}(n, \chi)$
is $B$-invariant and Zariski open in $X$.

\subsection{\it The quotient $X_B^{ss}(n,\chi)/\!/B$.}
\label{Bquotient} Denote the GIT quotient $$(X \times
G/B)_G^{ss}(L^n \otimes L_\chi)/\!/G$$ by $Q_{n,\chi}$ and let
$$\phi: (X \times G/B)_G^{ss}(L^n \otimes L_\chi) \rightarrow
Q_{n,\chi}$$ be the quotient map.  We then have the composition
map
$$\phi \circ \iota_B: X_B^{ss}(n,\chi)
\rightarrow (X \times G/B)_G^{ss}(L^n \otimes L_\chi) \rightarrow
Q_{n,\chi}.$$

\begin{thm}
\label{B-categoricalquotient} The morphism $\phi \circ \iota_B:
X_B^{ss}(n,\chi) \rightarrow Q_{n,\chi}$ is a categorical
quotient\footnote{For the definition of a categorical quotient,
see Definition 0.5 of \cite{GIT}.} for the $B$-action.
\end{thm}

\proof Let $\psi: X_B^{ss}(n,\chi) \rightarrow Z$ be any
$B$-morphism where $B$ acts trivially on $Z$. Then, one checks
that the map
$$\psi': (X \times G/B)_G^{ss}(L^n
\otimes L_\chi) \rightarrow Z$$
$$(x, g[B]) \rightarrow \psi(g^{-1} \cdot x)$$
is a $G$-morphism with respect to the trivial $G$-action on $Z$.
Clearly, $$\psi=\psi'\circ i_B.$$
 But
$$(X \times G/B)_G^{ss}(L^n
\otimes L_\chi)  \rightarrow Q_{n,\chi}$$ is categorical, hence we
have a commutative diagram
\begin{equation*}
\begin{CD}
  (X \times G/B)_G^{ss}(L^n \otimes L_\chi)  @>{\psi'}>> Z \\
  @V{\phi } VV    @V{\id} VV \\
   Q_{n,\chi} @>{\chi}>>  Z.
\end{CD}
\end{equation*}
This diagram extends to
\begin{equation*}
\begin{CD}
 X_B^{ss}(n,\chi) @>{i_B}>> (X \times G/B)_G^{ss}(L^n \otimes L_\chi)  @>{\psi'}>> Z \\
  && @V{\phi } VV    @V{\id} VV \\
  && Q_{n,\chi} @>{\chi}>>  Z
\end{CD}
\end{equation*}
which gives rise to the desired diagram
\begin{equation*}
\begin{CD}
  X_B^{ss}(n,\chi)  @>{\psi }>> Z \\
  @V{\phi \circ \iota_B } VV    @V{\id} VV \\
   Q_{n,\chi} @>{\chi}>>  Z.
\end{CD}
\end{equation*}
\endproof

Because of this theorem, we may also denote $Q_{n,\chi}$ by
$X_B^{ss}(n,\chi)/\!/B$.

\begin{lem}
\label{sectionIsom} (Guillemin-Sjamaar, \cite{GSj}) There is an
isomorphism of vector spaces
$$\rho: H^0(X \times G/B, L^{d} \otimes L_{d\chi})^G \rightarrow
H^0(X,L^d)^U_{d\chi}.$$
%which extends to an isomorphism of algebras
%$$\bigoplus_{d,\chi} H^0(X \times G/B, L^{d} \otimes L_{d\chi})^G \rightarrow
%\bigoplus_{, \chi} d H^0(X,L^d)^U_{d\chi}$$
\end{lem}
\proof
$$\rho: H^0(X \times G/B, L^{d} \otimes L_{d\chi})^G \rightarrow
H^0(X,L^d)^U_{d\chi}$$  is defined as follows. For any $\tilde{s}
\in H^0(X \times G/B, L^{d} \otimes L_{d\chi})^G$, then $s =
\rho(\tilde{s}) \in H^0(X, L^{d})^U$ is defined by
$$s(x) \otimes 1 = \tilde{s}(x,[B]), \forall x \in X.$$
One checks that so-defined $s$ is $U$-invariant and transforms
according to $d \chi$ under the action of the maximal torus $H$.
Conversely, a section $s \in H^0(X, L^{d})^U$ can be extended in a
unique way to a  section in $H^0(X \times G/B, L^{d} \otimes
L_{d\chi})^G$ by the formula
$$\tilde{s}(x, g[B]) = g(s(g^{-1} x) \otimes 1), \forall x \in X, g \in G.$$
\endproof

\begin{thm} We have
\label{MoreonBquotient}
%\begin{enumerate}
$$X_B^{ss}(n,\chi)=\{x \in X | \exists d > 0, s \in H^0(X, L^{nd})^{B[\chi]},
s(x) \ne 0\}.$$  In particular, the $B$-quotient
$X_B^{ss}(n,\chi)/\!/B$ is isomorphic to $\Proj
(R_{(n)}^{B[\chi]})$.
%\end{enumerate}
\end{thm}
\proof  By Lemma \ref{sectionIsom} (replace the line bundle $L$ by
$L^n$ in the lemma), we obtain an isomorphism
$$\rho: H^0(X \times G/B, L^{nd} \otimes L_{d\chi})^G \rightarrow H^0(X,L^{nd})^U_{d\chi}.$$
Since $(X \times G/B)_G^{ss}(L^n \otimes L_\chi)$ equals to
$$\{(x,g[B])| \exists d>0, \tilde{s} \in H^0(X \times G/B, L^{nd}
\otimes L_{d\chi})^G, \tilde{s}(x,g[B]) \ne 0\},$$ one checks from
the definition of $X_B^{ss}(n,\chi)$ that
$$X_B^{ss}(n,\chi)=\{x \in X | \exists d>0, s \in H^0(X,L^{nd})^U_{d\chi},
s(x) \ne 0\}.$$ Now observe that
$H^0(X,L^{nd})^U_{d\chi}=R^U_{{nd}, d\chi}$ is precisely the
subset of $B$-invariants of $R_{nd}$ under the $(-\chi)$-shifting,
that is,
\begin{equation}
\label{star} R^U_{nd, d\chi} = (R_{nd}^U)^{H[\chi]} =
R_{nd}^{B[\chi]} = H^0(X, L^{nd})^{B[\chi]}.
\end{equation}
This shows that
$$X_B^{ss}(n,\chi)=\{x \in X | \exists d>0, s \in H^0(X, L^{nd})^{B[\chi]},
s(x) \ne 0\}.$$

To show the last statement, note that the $B$-quotient
$X_B^{ss}(n,\chi) \rightarrow Q_{n,\chi}$ is identified with the
$G$-quotient
$$(X \times G/B)_G^{ss}(L^n \otimes L_\chi)/\!/G.$$ From the
above, we have that
$$\bigoplus_d H^0(X \times G/B, L^{nd} \otimes L_{d\chi})^G =
\bigoplus_d R_{nd}^{B[\chi]} = R_{(n)}^{B[\chi]}.$$
Because $(X \times G/B)_G^{ss}(L^n \otimes L_\chi)/\!/G$ is
isomorphic to
$$\Proj (\bigoplus_d H^0(X \times G/B, L^{nd} \otimes L_{d\chi})^G
),$$ we obtain that the $B$-quotient $Q_{n,\chi}$ is isomorphic to
$\Proj (R_{(n)}^{B[\chi]})$.
\endproof

We isolate the following identity from Equation (\ref{star}) in
the
proof of the above theorem. %Theorem \ref{MoreonBquotient}.

\begin{cor}
\label{isolatedid}
 $R_{(n)}^{B[\chi]} = (R_{(n)}^U)^{H[\chi]}$.
\end{cor}

\section{Proof of Theorem \ref{mainthm}}

For ${\chi \over n} \in \Delta_{\rm reg}$, the existences of the
three quotients
$$(X^{ss}_U(L)/\!/U)_H^{ss}(n,\chi)/\!/H,$$ $$(X \times G/B)_G^{ss}(L^n
\otimes L_\chi)/\!/G,$$ $$X_B^{ss}(n, \chi)/\!/B,$$ are proved in
\S\S \ref{torusaction}, \ref{byG}, \ref{byB}, respectively.

That $X_B^{ss}(n, \chi)/\!/B$ is isomorphic to $(X \times
G/B)_G^{ss}(L^n \otimes L_\chi)/\!/G$ is contained in Theorem
\ref{B-categoricalquotient}.

To show that $X_B^{ss}(n, \chi)/\!/B$ is isomorphic to
$(X^{ss}_U(L)/\!/U)_H^{ss}(n,\chi)/\!/H,$ note that by Theorem
\ref{MoreonBquotient}, $X_B^{ss}(n, \chi)/\!/B$ is isomorphic to
$\Proj (R_{(n)}^{B[\chi]})$. By Corollary \ref{isolatedid}, it is
isomorphic to $\Proj ((R_{(n)}^U)^{H[\chi]}). $ Now it follows
from Theorem \ref{Hquotient} that it is isomorphic to
$(X^{ss}_U(L)/\!/U)_H^{ss}(n,\chi)/\!/H.$

\section{Singular rational points of $\Delta$ and $G/P$.}
\label{G/P}

\subsection{\it The action $G \times (X \times G/P) \rightarrow (X \times
G/P)$} For a rational point ${\chi \over n} \in \Delta$ that lies
on a wall of the Weyl chamber, the character $\chi \in \Lambda_+$
determines a parabolic subgroup $P$ strictly containing $B$:
$$P=\{ g \in G | \lim_{t \to 0} \chi (t) g \chi (t)^{-1} \; {\rm
exists} \}.$$
%the line bundle $L_\chi = G \times_B \C_{-\chi}$ is not ample. For
%a suitable parabolic subgroup $P$ strictly containing $B$, if we
Let $\C_{-\chi}$ be the one-dimensional $P$-module with the
character $-\chi$. Then $L_\chi' = G \times_P \C_{-\chi}$ is a
$G$-linearized ample line bundle over $G/P$.
% and $L_\chi$ is the pull-back of $L_\chi'$ by the projection $\pi: G/B
%\rightarrow G/P$.

%For the linearized line bundle $L^n \otimes L_\chi $, the
%$G$-semistable locus on $X \times G/B$ is empty, so the
%$G$-quotient on $X \times G/B$ is empty. But the algebra
%$R_{(n)}^{B[\chi]}$ may well be nonzero, that is, the $B$-quotient
%$\Proj(R_{(n)}^{B[\chi]})$ may still be a nonempty variety.
To extend the correspondences of Theorem \ref{mainthm} to this
case, we can simply replace the second action by the diagonal
action
$$G \times (X \times G/P) \rightarrow (X \times G/P).$$

\subsection{\it Extensions of some results of \S \ref{byB}.}
Lemma \ref{sectionIsom}, with basically the same proof
(\cite{GSj}), now reads: we have an isomorphism of vector spaces
\begin{equation}
\label{eqForG/P} H^0(X \times G/P, L^{d} \otimes L_{d\chi}')^G
\rightarrow H^0(X,L^d)^U_{d\chi}
\end{equation}
where $L_{d\chi}'=(L_\chi')^d$.

 Equation (2) in the proof of Theorem
\ref{MoreonBquotient} remains true without any change.
%hence, so does Corollary \ref{isolatedid}.

Since the $G$-quotient $(X \times G/P)_G^{ss}(L^n \otimes
 L_\chi')/\!/G$ is isomorphic to
$$\Proj (\bigoplus_d H^0(X \times G/P, L^{nd} \otimes L'_{d\chi})^G
),$$  by Equation (\ref{eqForG/P}) of this section and Equation
(2) in the proof of Theorem \ref{MoreonBquotient}, we will obtain
the following.

%Analogous assertions in \S\S \ref{B-opensubset} and
%\ref{Bquotient} remain true when $B$, $G/B$ and $L_\chi$ are
%(carefully) replaced by $P$, $G/P$ and $L_\chi'$.  Then, we will
%obtain that
%we can take the following approach, utilizing the Relative
%Geometric Invariant Theory (\cite{rgit}). Assume that $(X \times
%G/P)^{ss}(L_\chi') = (X \times G/P)^s (L_\chi')$, then for a {\it
%small perturbation} $L_{\chi, {\rm small}}$ of the linearization
%$L_\chi = \pi^* L_\chi'$, we have that $$(X \times
%G/B)^{ss}(L_{\chi, {\rm small}}) = (X \times G/B)^s(L_{\chi, {\rm
%small}})= \pi^{-1}( (X \times G/P)^s (L_\chi'))$$ (see \cite{rgit}
%for details).

\subsection{\it The correspondences.}
\begin{thm}
\label{thm:G/P}
The $G$-quotient $(X \times G/P)_G^{ss}(L^n
\otimes
 L_\chi')/\!/G$,
the $B$-quotient $\Proj(R_{(n)}^{B[\chi]})$,  and  the
$H$-quotient $\Proj((R_{(n)}^U)^{H[\chi]})$ are isomorphic to each
other.
\end{thm}

\begin{rem}
\label{chi=0}
 The case when $\chi = 0$ is worth mentioning. In this case, the
 parabolic subgroup is $G$ so that $G/P$ is a point,
 hence $(X \times G/P)_G^{ss}(L^n \otimes
 L_\chi')/\!/G$ is just  the $G$-quotient $X^{ss}_G (L)/\!/G$.
 The fact that the $G$-quotient $X^{ss}_G (L)/\!/G$, the $B$-quotient $\Proj(R_{(n)}^{B[0]})$
 and the $H$-quotient $\Proj((R_{(n)}^U)^{H[0]})$ are all isomorphic
 can also be seen by observing that
$$R_{(n)}^G = R_{(n)}^{B[0]}= (R_{(n)}^U)^{H[0]}.$$
\end{rem}

Replacing $B$ by $P$ in \S \ref{B-opensubset}, we will obtain a
$P$-invariant Zariski open subset $X^{ss}_P(n, \chi)$ of $X$. Then
a proof almost exactly the same as that of Theorem
\ref{B-categoricalquotient} will yield the following (details are
left to the reader).

\begin{thm}
\label{P-categoricalquotient} The morphism
$$ X_P^{ss}(n,\chi) \rightarrow (X \times G/P)_G^{ss}(L^n \otimes
 L_\chi') \rightarrow (X \times G/P)_G^{ss}(L^n \otimes
 L_\chi')/\!/G$$ is a categorical quotient for the $P$-action on $X$.
\end{thm}

\section{\it Concluding remarks}
\label{concluding}

\subsection{\it Singularities of the unipotent quotient.}
We know little about the singularities of the unipotent quotient
$X^{ss}_U(L)/\!/U$. However, the correspondences  of Theorem
\ref{mainthm} shed some lights on it.

When ${\chi \over n } \in \Delta_{\rm reg}$ is a regular value of
the moment map $\Phi$, by Remark \ref{orbifolds}, the $G$-quotient
$$(X \times G/B)_G^{ss}(L^n \otimes L_\chi)/\!/G$$ is an orbifold,
that is, it has at worst finite quotient singularities. By Theorem
\ref{mainthm}, the same hold for the corresponding $H$-quotient
$$(X^{ss}_U(L)/\!/U)_H^{ss}(n,\chi)/\!/H.$$ This indicates that the
Zariski open subset $(X^{ss}_U(L)/\!/U)^{ss}(n,\chi)$ of the
unipotent quotient $X^{ss}_U(L)/\!/U$ has at worst finite quotient
singularities, and this holds for all almost all rational points
${\chi \over n} \in \Delta$. The variety $X^{ss}_U(L)/\!/U$ and
its applications call for further investigation.

\subsection{\it More on Example \ref{1stexmp}.}
\label{moreonexample}
 For the line bundle
 $L = \otimes_i \cO(d_i)$ over $(\proj^1)^n$ with $d_n$ sufficiently
 large relative to other $d_i$ ($1 \le i \le n-1$),
we have that the unipotent quotient is isomorphic to
$(\proj^1)^{n-1}.$ Note that in this case the homogeneous space
$G/B$ is isomorphic to $\proj^1$. One checks that the three
actions in this case are the following diagonal actions
$$H \times (\proj^1)^{n-1} \rightarrow (\proj^1)^{n-1},$$
$$G \times (\proj^1)^{n+1} \rightarrow (\proj^1)^{n+1},$$
$$B \times (\proj^1)^n \rightarrow (\proj^1)^n.$$
For any ${\chi \over m} \in \Delta_{\rm reg}$, the corresponding
quotient of the first action is a toric variety\footnote{GIT
quotients of a projective toric variety by a subtorus are again
projective toric varieties (\cite{KSZ}).}, hence so is the
corresponding quotient of the second action by Theorem
\ref{mainthm}. This implies that the $G$-linearized line bundle
$L^m \otimes L_\chi$ over $(\proj^1)^{n+1}$, with $d_n$
sufficiently
 large relative to the rest, is a very special one, because
for a general ample $G$-linearized line bundle over
$(\proj^1)^{n+1}$ ($n \ge 4$), we know that the corresponding GIT
quotient is not toric. For example, when $n=4$, with respect to
the $G$-linearized line bundle $\otimes_{i=1}^5 \cO(1)$, the GIT
quotient of $(\proj^1)^{4+1}$ is isomorphic to the blowup of
$\proj^2$ along 4 general points which is not toric. By our main
correspondences, this implies that for the line bundle
$\otimes_{i=1}^n \cO(d_i)$ with general positive integers $d_i$
($1 \le i \le n$), the corresponding unipotent quotient of
$(\proj^1)^n$ can not be toric variety. In particular, it is not
isomorphic to the unipotent quotient $(\proj^1)^{n-1}.$ This
justifies the assertion of Remark \ref{distintqUq} that the
unipotent quotient $X^{ss}_U(L)/\!/U$, in general, depends on the
choice of the underlying line bundle $L$.

It is an interesting problem to (explicitly) determine
$X^{ss}_U(L)$ and $X^{ss}_U(L)/\!/U$ for general choices of $d_i
(1 \le i \le n)$.

Finally, we mention that the  GIT quotients of the second action
here can be interpreted as moduli spaces of spacial polygons
(\cite{KM}). We do not know whether the other two admit natural
geometric explanations.

\subsection{\it Related and further works.}
There are a number of papers (e.g., \cite{F}, \cite{GP1},
\cite{GP2}, \cite{Snow}) that study quotients of unipotent group
actions or quotients of general algebraic group actions to which
this paper is related.

There are some moduli spaces that may be constructed as quotients
of Borel subgroups. For example, the moduli spaces of vector
bundles over smooth algebraic curves with {\it complete} parabolic
structures are naturally quotients by Borel subgroups (see page
545 of \cite{BH}. For partial  parabolic structures, one should
use parabolic subgroups instead). Via a shifting trick similar to
that of \S \ref{shiftingtrick}, these moduli spaces  are
constructed as quotients by reductive groups by Mehta and Sashadri
in \cite{MS}. Our work here indicates that they may also be
constructed as quotients by torus actions. This would use certain
{\it unipotent quotients}. Thus it would be an interesting problem
to see what these unipotent quotients are and whether they admit
interesting moduli interpretations.

When $X_G^{ss}(L)=X^s_G(L)$ (cf. Remark \ref{chi=0}), Brion
proposed the following: through (orbifold) fiber bundle and toric
flips, we may relate the quotient $X_G^{ss}(L)/\!/G$ by the
reductive group $G$ to a quotient of $X^{ss}_U(L)/\!/U$ by the
maximal torus $H$. This would give an alternative way to study the
topology of a general GIT quotient (cf. \cite{Kirwan84}). Some
related works around this area may be found in \cite{AC} and
\cite{BC}.

After receiving the preliminary version of this paper, Brion
mentioned to me that he was also convinced that the results here
hold. Part of his idea appeared in L. Pillons' thesis
\cite{Pillons}.

\bigskip

{\sl Acknowledgements.} Many works have inspired this article,
including Mumford's appendix to Ness's paper (\cite{Ness}),
Brion's paper on moment map images (\cite{Brion}), and many
others' works on symplectic reductions through general co-adjoint
orbits. The author, very possibly among many others in the area,
was convinced that there are natural correspondences and
isomorphisms among the quotients of the three actions quite a
while ago. Finally, Guillemin-Sjamaar's article (\cite{GSj})
convinced him that it should be written up.  I thank Michel Brion
for very helpful correspondences, comments and corrections. I am
also very grateful to the referee for providing many corrections
and very critical, useful comments. I am especially thankful to
the editors for calling my attention to the quality of the
exposition. As a result, the paper was completely re-organized and
the bulk of it was more carefully re-written.

\end{document}